# HIGHER CRITICISM FOR DETECTING SPARSE HETEROGENEOUS MIXTURES[1]

BY DAVID DONOHO AND JIASHUN JIN

*Stanford University*

*In Memory of John Wilder Tukey, 1915–2000*

*Higher criticism*, or *second-level significance testing*, is a multiple-comparisons concept mentioned in passing by Tukey. It concerns a situation where there are many independent tests of significance and one is interested in rejecting the joint null hypothesis. Tukey suggested comparing the fraction of observed significances at a given $\alpha$-level to the expected fraction under the joint null. In fact, he suggested standardizing the difference of the two quantities and forming a $z$-score; the resulting $z$-score tests the significance *of the body of significance tests*.

We consider a generalization, where we maximize this $z$-score over a range of significance levels $0 < \alpha \leq \alpha_0$. We are able to show that the resulting *higher criticism statistic* is effective at resolving a very subtle testing problem: testing whether $n$ normal means are all zero versus the alternative that a small fraction is nonzero.

The subtlety of this "sparse normal means" testing problem can be seen from work of Ingster and Jin, who studied such problems in great detail. In their studies, they identified an interesting range of cases where the small fraction of nonzero means is so small that the alternative hypothesis exhibits little noticeable effect on the distribution of the $p$-values either for the bulk of the tests or for the few most highly significant tests. In this range, when the amplitude of nonzero means is calibrated with the fraction of nonzero means, the likelihood ratio test for a precisely specified alternative would still succeed in separating the two hypotheses.

We show that the higher criticism is successful throughout the same region of amplitude sparsity where the likelihood ratio test

Received August 2002; revised May 2003.

[1]Supported in part by National Science Foundation Grants DMS 00-77261 and DMS 95-05151.

would succeed. Since it does not require a specification of the alternative, this shows that higher criticism is in a sense optimally adaptive to unknown sparsity and size of the nonnull effects. While our theoretical work is largely asymptotic, we provide simulations in finite samples and suggest some possible applications. We also show that higher critcism works well over a range of non-Gaussian cases.

**1. Introduction.** In his Class Notes for Statistics 411 at Princeton University in 1976 [31], Tukey introduced the notion of the *higher criticism* by means of a story. A young psychologist administers many hypothesis tests as part of a research project, and finds that, of 250 tests 11 were significant at the 5% level. The young researcher feels very proud of this fact and is ready to make a big deal about it, until a senior researcher (Tukey himself?) suggests that one would expect 12.5 significant tests even in the purely null case, merely by chance. In that sense, finding only 11 significant results is actually somewhat disappointing!

Tukey used this story as a way to make vivid the notion of the *higher criticism* of such situations as multiple testing. He then proposed a sort of *second-level significance testing*, based on the statistic

$$\mathrm{HC}_{0.05,n} = \sqrt{n}[(\text{Fraction Significant at } 0.05) - 0.05]/\sqrt{0.05 \times 0.95},$$

and suggested that values of (say) 2 or greater indicate a kind of *significance of the overall body of tests*. (The same statistic was proposed and applied in a psychometric trial by Brožek and Tiede [10] even earlier, but without the catchy name.)

Although Tukey's discussion turned to other topics at that point, we may, if we like, imagine that it had continued in this vein. We might then consider not only significance at the 0.05 level, but perhaps at all levels between (say) 0 and $\alpha_0$, and so define

$$\mathrm{HC}_n^* = \max_{0 < \alpha \leq \alpha_0} \sqrt{n}\,[(\text{Fraction Significant at } \alpha) - \alpha]/\sqrt{\alpha \times (1-\alpha)}.$$

In this paper, we will analyze a statistic of this kind in a setting where there are a small fraction of nonnull hypotheses and derive an adaptive optimality for it.

In our setting there are $n$ independent tests of unrelated hypotheses, $H_{0,i}$ vs. $H_{1,i}$, where the test statistics $X_i$ obey

$$(1.1) \qquad H_{0,i} : X_i \sim N(0,1),$$

$$(1.2) \qquad H_{1,i} : X_i \sim N(\mu_i, 1), \qquad \mu_i > 0.$$

In the overwhelming majority of the tests, the corresponding null hypothesis is true (i.e., the corresponding normal mean $\mu_i = 0$), but some small fraction *may* concern tests where the null hypothesis is false (and so the mean $\mu_i > 0$).

HIGHER CRITICISM FOR DETECTING MIXTURES 3So the fraction of false null hypotheses, if nonzero, is small. Can we tell reliably whether the fraction is actually 0 or not?

We mention three application areas where situations like this might arise:

- *Early detection of bioweapons use.* Suppose that there are $n$ observational units in a certain geographical region, and for each one we have a $z$-score associated with the presence of a certain symptom at rates higher than background. If a bioweapon has been used in that region, then in early stages we do not expect all observational units to be affected, we do not know which ones might be affected, and we do not want to wait until some observational unit begins to display wildly elevated rates. We want to detect while a small fraction begins to show individually significant results but no unit yet shows jointly significant results.
- *Detection of covert communications.* In a signals intelligence setting we suppose that a small fraction of the signal spectrum in a certain situation may be used for covert communications, which would mean that a few frequencies exhibit increased power. However, we do not know what frequencies those might be, and the specific frequencies being used might change randomly from one epoch to another, so that we never get a very definite indication that we are definitely seeing increased power in any one specific frequency. Nevertheless, we might still want to detect the presence of a small fraction of frequencies with slightly increased power.
- *Meta-analysis with heterogeneity.* We have results from $n$ experiments testing a certain treatment. It turns out that an unidentified experimental factor is crucial to success, but that only in a small fraction of experiments is this factor fortuitously chosen so that the experimental performance follows a nonnull distribution. Can we reliably detect the presence of a small fraction of well-laid out experiments among many hopeless ones, when we do not know which ones may be well-laid out?

There are many other potential applications in signal processing; see, for example, [19]–[21]. In spatial statistics Kendall and Kendall [25] developed a statistic closely related to $\mathrm{HC}_n^*$ called "pontogram" for the purpose of detecting near-alignments in sets of points.

1.1. *The model, and the asymptotic detection boundary.* Translating our problem into precise terms, we begin by scrutinizing a special case where all the nonzero $\mu_i$ are equal, and we can then model our data as providing $n$ i.i.d. observations from one of two possible situations:

(1.3) $\qquad H_0 : X_i \stackrel{\text{i.i.d.}}{\sim} N(0,1), \qquad 1 \leq i \leq n,$

(1.4) $\qquad H_1^{(n)} : X_i \stackrel{\text{i.i.d.}}{\sim} (1-\varepsilon)N(0,1) + \varepsilon N(\mu,1), \qquad 1 \leq i \leq n.$

Here $H_0$ denotes the global intersection null hypothesis, and $H_1^{(n)}$ denotes a specific element in its complement. Under $H_1^{(n)}$, a fraction $\varepsilon$ of the data comes from a normal with common nonnull mean. Here $\varepsilon = \varepsilon_n$ and the mean $\mu = \mu_n$ will be chosen to make the problem very hard, but (just barely) still solvable.

Obviously, in this situation, with $\varepsilon_n$ and $\mu_n$ fixed and known, the optimal procedure is simply the likelihood ratio test; a careful analysis of its performance in [23] (cf. also [16] and [17]) tells us the following. Suppose we let $\varepsilon_n = n^{-\beta}$ for some exponent $\beta \in (\frac{1}{2}, 1)$, so that the fraction of nonzero means is small but not vanishingly small. In this range the number of nonzero means is too small to be noticeable in any sum which is, in expectation, of order $n$, so it cannot noticeably affect the behavior of the bulk distribution of the $p$-values. Let

$$(1.5) \qquad \mu_n = \sqrt{2r \log(n)}, \qquad 0 < r < 1.$$

As $\mu_n < \sqrt{2 \log(n)}$, the nonzero means are, in expectation, smaller than the largest $X_i$ coming from the true component null hypotheses, so the nonzero means cannot have a visible effect on the upper extremes. Clearly, this is a rather subtle testing problem.

It turns out that there is a *threshold effect* for the likelihood ratio test: the sum of Type I and Type II errors tends to 0 or 1 depending on whether $\mu$ exceeds a so-called *detection boundary* or not. In detail, there is a function $\rho^*(\beta)$ so that

$$\text{if } r > \rho^*(\beta), \qquad H_0 \text{ and } H_1^{(n)} \text{ separate asymptotically,}$$
$$\text{if } r < \rho^*(\beta), \qquad H_0 \text{ and } H_1^{(n)} \text{ merge asymptotically.}$$

In short, $\rho^*(\beta)$ defines a precise demarcation between what is possible and impossible in this problem, that is, how big the nonzero effect must be to be detectable as a function of the rarity of nonzero effects. Hence, we have the term *detection boundary*. Indeed, translating results of Ingster [17] to our notation (see also [23]),

$$(1.6) \qquad \rho^*(\beta) = \begin{cases} \beta - \frac{1}{2}, & \frac{1}{2} < \beta \leq \frac{3}{4}, \\ (1 - \sqrt{1-\beta})^2, & \frac{3}{4} < \beta < 1. \end{cases}$$

If we think of the $(r, \beta)$ plane, $0 < r < 1$, $\frac{1}{2} < \beta < 1$, we are saying that, throughout the region $r > \rho^*(\beta)$ the alternative can be detected reliably using the likelihood ratio test (LRT). Unfortunately, the usual (Neyman–Pearson) likelihood ratio requires a precise specification of $r$ and $\beta$, and misspecification of $(r, \beta)$ may lead to failure of the LRT; see [23] for a discussion. Naturally, in any practical situation we would like to have a procedure which does well throughout this whole region without knowledge of $(r, \beta)$.

Bickel and Chernoff [6] and Hartigan [14] have shown that the usual generalized likelihood ratio test $\max_{\varepsilon,\mu}(dP_1^{(n)}(\varepsilon,\mu)/dP_0^{(n)})(X)$ has nonstandard behavior in this setting; in fact the maximized ratio tends to $\infty$ under $H_0$. It is not clear that this test can be relied on to detect subtle departures from $H_0$. Ingster [18] has proposed an alternative method of adaptive detection which maximizes the likelihood ratio over a finite but growing list of simple alternative hypotheses. By careful asymptotic analysis, he has, in principle, completely solved the problem of adaptive detection in this setting; however, this is a relatively complex and delicate procedure which is tightly tied to the narrowly specified model (1.3) and (1.4). It would be nice to have an easily implemented and intuitive method of detection which is able to work effectively throughout the whole region $0 < \beta < \frac{1}{2}$, $r > \rho^*(\beta)$, which is not tied to the narrow model (1.3) and (1.4), and which is in some sense easily adapted to other (non-Gaussian) mixture models. This is where $\mathrm{HC}_n^*$ comes in.

1.2. *Performance of higher criticism.* To apply the higher criticism, let us convert the individual test statistics into another form. Let $p_i = P\{N(0,1) > X_i\}$ be the $p$-value for the $i$th component null hypothesis, and let the $p_{(i)}$ denote the $p$-values *sorted in increasing order*, so that under the intersection null hypothesis the $p_{(i)}$ behave like order statistics from a uniform distribution.

With this notation, we can write

$$\mathrm{HC}_n^* = \max_{1 \leq i \leq \alpha_0 \cdot n} \sqrt{n}\,[i/n - p_{(i)}]/\sqrt{p_{(i)}(1 - p_{(i)})}.$$

Despite the closeness of our statistic to the one in [2], note that what we are doing is not arbitrary goodness of fit; instead we are dealing with a specific kind of multiple hypothesis testing, outside of which the problem would not be interesting.

To use $\mathrm{HC}_n^*$ to conduct a level-$\alpha$ test, we must find a critical value $h(n, \alpha)$:

$$P_{H_0}\{\mathrm{HC}_n^* > h(n, \alpha)\} \leq \alpha.$$

Adapting asymptotic theory for the normalized empirical process as in [27], Chapter 16, gives us the following information on the size of $h(n, \alpha)$:

THEOREM 1.1. *Under the null hypothesis $H_0$,*

(1.7) $$\frac{\mathrm{HC}_n^*}{\sqrt{2 \log \log(n)}} \xrightarrow{p} 1, \qquad n \to \infty.$$

It follows that, for fixed $\alpha > 0$, $h(n, \alpha) \approx \sqrt{2 \log \log(n)}$. For asymptotic analysis, it is convenient to consider a sequence of problems indexed by $n$, with critical values $\alpha_n \to 0$. We will say that the level $\alpha_n \to 0$ *slowly enough* if

$$h(n, \alpha_n) = \sqrt{2 \log \log(n)}(1 + o(1)).$$

THEOREM 1.2. *Consider the higher criticism test that rejects $H_0$ when*

$$\mathrm{HC}_n^* > h(n, \alpha_n),$$

*where the level $\alpha_n \to 0$ slowly enough. For every alternative $H_1^{(n)}$ defined above where $r$ exceeds the detection boundary $\rho^*(\beta)$—so that the likelihood ratio test would have full power—higher criticism also has full power:*

$$P_{H_1^{(n)}}\{\text{Reject } H_0\} \to 1, \qquad n \to \infty.$$

Roughly speaking, everywhere in the amplitude/sparsity $(r, \beta)$ plane—where the likelihood ratio test would completely separate the two hypotheses asymptotically—the higher criticism will also completely separate the two hypotheses asymptotically. Of course, in the cases where the amplitude/sparsity relation falls *below* the detection boundary, all methods fail. More precisely, we only claim that higher criticism works in the interior of this region. Just at the critical point where $r = \rho^*(1 + o(1))$, our result says nothing; this would be an interesting (but very challenging) area for future work.

1.3. *Which part of the sample contains the information?* Underlying our results is a set of insights about "where to look" for evidence against $H_0$ [30], how the evidence may not be in the "obvious" place and how the adaptation in HC* automatically ensures that the best evidence will be included in making the "case" against $H_0$.

To get started, note that, in the null case $\mathrm{HC}_n^*$ is closely related to well-known functionals of the standard uniform empirical process. This is because, under $H_0$, the $n$ $p$-values can be viewed as i.i.d. samples from $U(0, 1)$.

Formalizing, given $n$ independent random samples $U_1, \ldots, U_n$ from $U(0, 1)$, with empirical distribution function

$$F_n(t) = \frac{1}{n} \sum_{i=1}^n \mathbb{1}_{\{U_i \leq t\}},$$

the uniform empirical process is denoted by

$$U_n(t) = \sqrt{n}\,[F_n(t) - t], \qquad 0 < t < 1,$$

and the *normalized uniform empirical process* by

$$W_n(t) = \frac{U_n(t)}{\sqrt{t(1-t)}}.$$

Note that $W_n(t)$ is asymptotically $N(0,1)$ for each fixed $t \in (0,1)$.

Since, under $H_0$ the $p$-values are i.i.d. $U(0,1)$, we have, in that case, the representation

$$\mathrm{HC}_n^* = \max_{0 < t < \alpha_0} W_n(t);$$

note the use of max instead of sup, since the maximum value of $W_n(t)$ is attained.

We will extend usage below by letting $W_n(t)$ also stand for the normalized empirical process starting from $F_n(t) = \frac{1}{n} \sum_{i=1}^n \mathbb{1}_{\{p_i \le t\}}$, where the $p_i$ are the $n$ given $p$-values, but which are i.i.d. $U(0,1)$ only in the null case. Accordingly, $W_n(t)$ will be asymptotically $N(0,1)$ at each fixed $t$ only under $H_0$, and one anticipates a different limiting behavior under the alternative.

We now look for the value of $t$ at which $W_n(t)$ has the most dramatically different behavior under the null and the alternative.

We introduce the notation $z_n(q) = \sqrt{2q \log(n)}$, for $0 < q \le 1$, and look for the value of $q$ for which $\mathrm{Prob}\{X_i > z_n(q)\}$ best differentiates between null and alternative. Recall that, under the alternative $H_1^{(n)}$ the data have a sparse component with nonzero mean at $\mu_n = z_n(r)$. It might seem that the most informative part of the sample would be in the vicinity of $\mu_n$, where data from the alternative are most common, and that therefore the most informative value of $q$ is $q = r$. Surprisingly, this is not the case.

To find this most informative $q$, we introduce some notation. Let $p_{n,q} = P\{N(0,1) > z_n(q)\}$, $q > 0$, and note that

$$\max_{0 < t < 1/2} W_n(t) = \max_{0 < q < \infty} W_n(p_{n,q}).$$

It will be immediately clear that it is only necessary to consider $0 < q \le 1$.

Let $N_n(q)$ count the observations exceeding $z_n(q)$:

$$N_n(q) = \#\{i : X_i \ge z_n(q)\}.$$

Then define

$$V_n(q) = \frac{N_n(q) - n p_{n,q}}{\sqrt{n p_{n,q}(1 - p_{n,q})}}.$$

In many calculations, we will find factors with polylog behavior, $\sim \mathrm{Const} \log^a(n)$, where $a$ may be positive or negative depending on the case. When such factors are multiplying terms $n^\gamma$ for $\gamma \ne 0$, the polylog factors have a weak influence on the eventual growth rate of the resulting product. To focus on the main ideas, we introduce the notation $L_n$ for a generic polylog factor, which may change from occurrence to occurrence. When we do this, we think of $L_n$ essentially as if it were constant.

Now

$$P\{N(\mu_n, 1) > z_n(q)\} = L_n \cdot n^{-(\sqrt{q}-\sqrt{r})^2},$$
$$P\{N(0, 1) > z_n(q)\} = L_n \cdot n^{-q}, \qquad r < q \leq 1.$$

It follows that, under the alternative $H_1^{(n)}$, we have

$$EV_n(q) = L_n \cdot \frac{n^{1-\beta} \cdot n^{-(\sqrt{q}-\sqrt{r})^2}}{\sqrt{n \cdot n^{-q}}} = L_n \cdot n^{[(1+q)/2-\beta-(\sqrt{q}-\sqrt{r})^2]},$$

while under the null $EV_n(q) = 0$. The most informative value of $q$ will optimize the growth rate of $n^{[(1+q)/2-\beta-(\sqrt{q}-\sqrt{r})^2]}$ to $\infty$. The most informative value of $q$ satisfies

$$\text{if } r < \tfrac{1}{4}, \qquad \text{then } q = 4r \text{ and } EV_n(q) = L_n \cdot n^{[r-(\beta-1/2)]},$$
$$\text{if } r \geq \tfrac{1}{4}, \qquad \text{then } q = 1 \text{ and } EV_n(q) = L_n \cdot n^{[(1-\beta)-(1-\sqrt{r})^2]}.$$

In fact, there can be considerable latitude in choosing $q$ so that $EV_n(q)$ goes to $\infty$ under $H_1^{(n)}$. This requires

$$2\sqrt{r} - \sqrt{2(r-\beta+\tfrac{1}{2})} < \sqrt{q} < 2\sqrt{r} + \sqrt{2(r-\beta+\tfrac{1}{2})}.$$

Notice that $r > \rho^*(\beta)$ implies $r - \beta + \tfrac{1}{2} > 0$. Within this interval the "most informative choice" of $q$ is the center: $\sqrt{q} = 2\sqrt{r}$ in case $r < \tfrac{1}{4}$.

Translating this into the original $z$-scale yields the surprise mentioned earlier. When $r \leq \tfrac{1}{4}$, the most informative place on the original $z$-scale is not at $\mu_n$, as we might suppose, but at $2\mu_n$. By going out "in the tails" farther than $\mu_n$, observations from both $H_0$ and $H_1^{(n)}$ are becoming extremely rare, but the ones from $H_1^{(n)}$ are far more frequent in a relative sense.

The story when $r > \tfrac{1}{4}$ is somewhat different. There, when looking for a discrepancy, we had best look near $\sqrt{2\log(n)}$, which is less than $2\mu_n$. The point is that observations from $H_0$ almost never get substantially larger than $\sqrt{2\log(n)}$, so there is no need to look much farther out in the tails.

1.4. *Comparison to several multiple comparison procedures.* The higher criticism is just one specific approach to combining many $p$-values in search of an overall test; many other tools are available from the field of multiple comparisons (e.g., [15], [26] and [33]) and meta-analysis [3]. How do these other tools perform?

In this section we describe several specific procedures and their range of detectability. See Figure 1.

1.4.1. *Range/Maximum/Bonferroni.* One of the most classical and frequently used tools in multiple comparisons, also associated with Tukey [33], is the studentized range:

$$R_n = (\max X_i - \min X_i)/S_n,$$

where $S_n$ is the sample standard deviation. This is frequently used in testing for homogeneity of a set of normal means, and could well be used in the setting we have defined here. For our theoretical purposes, it is convenient to analyze the simpler statistic

$$M_n = \max_i X_i,$$

where we focus attention on one-sided deviations only and have no need to estimate the (known) standard deviation under the null. (We note that in the field of meta-analysis this is the same thing as combining several $p$-values by taking the minimum $p$-value [3]; another equivalent term for this is Bonferroni-based inference.) The maximum statistic has a critical value $m(n, \alpha)$ which obeys

$$m(n, \alpha) \sim \sqrt{2 \log(n)}, \qquad n \to \infty.$$

In comparison to HC$^*$, it follows that the test focuses entirely on whether there are observations exceeding $\sqrt{2 \log(n)}$. We again say that $\alpha_n$ goes to 0 *slowly enough* (now for use with $M_n$) if

$$m(n, \alpha_n) \sim \sqrt{2 \log(n)}.$$

The following result summarizes the behavior of $M_n$.

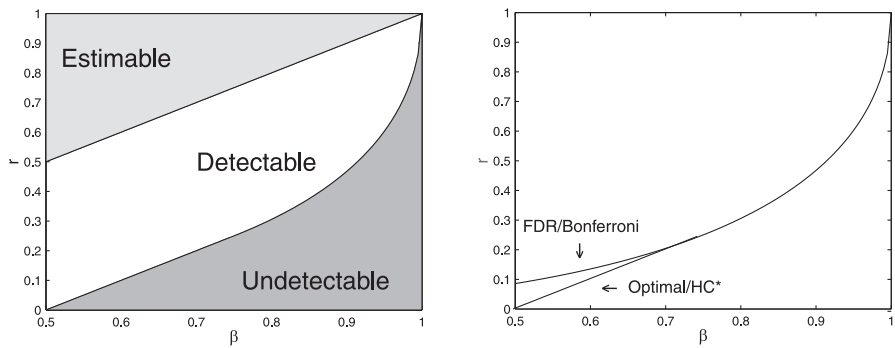

FIG. 1. Left: *three regions of the $\beta - r$ plane. The detection boundary separates the detectable region from the undetectable region. For the estimable region, it is possible not only to detect the presence of nonzero means, but also to estimate those means.* Right: *two detection boundaries. The one on the bottom is the optimal detection boundary as well as the detection boundary for* HC$^*$; *the one on the top is for range/maximum/Bonferroni/FDR. Two detection boundaries are only different when $\frac{3}{4} < \beta < 1$.*

THEOREM 1.3. *Define*

$$\rho_{\mathrm{Max}}(\beta) = (1 - \sqrt{1-\beta})^2.$$

*Suppose $r > \rho_{\mathrm{Max}}(\beta)$ and consider a sequence of level $\alpha_n$-tests based on $M_n$ with $\alpha_n \to 0$ slowly enough. Then*

$$P_{H_1^{(n)}}\{M_n \text{ rejects } H_0\} \to 1.$$

In short, $\rho_{\mathrm{Max}}$ defines the detection boundary for $M_n$. This compares to the "efficient" boundary as follows:

$$\rho^*(\beta) = \rho_{\mathrm{Max}}(\beta), \qquad \beta \in [\tfrac{3}{4}, 1),$$

so that $M_n$ is effective in the range of very sparse alternatives, while

$$\rho^*(\beta) < \rho_{\mathrm{Max}}(\beta), \qquad \beta \in [\tfrac{1}{2}, \tfrac{3}{4}),$$

so that $M_n$ is inefficient if $\beta < \tfrac{3}{4}$. In particular, note that

$$0 = \rho^*(\tfrac{1}{2}) < \rho_{\mathrm{Max}}(\tfrac{1}{2}) = (2 - \sqrt{2})^2/4 \approx 0.0858.$$

We interpret this as saying that $M_n$ (and the Studentized range $R_n$ as well) can be dramatically outperformed when there may be about $n^{1/2}$ nonzero means among the $n$ observations. Compare Figure 1.

1.4.2. *FDR-controlling methods.* Recently considerable interest has been focused on the false discovery rate (FDR)-controlling methodology for simultaneous inference. In one example of this so-called FDR approach [4], one considers, for $k = 1, 2, 3, \ldots$, the $k$ most significant $p$-values. These are compared with $\alpha \frac{k}{n}$, where $\alpha$ is a critical value (e.g., 0.05). If the $p$-values in the group are all smaller than the standard for comparison, then $H_0$ is accepted at that stage. If some are larger than the standard for comparison, then $H_1^{(n)}$ is accepted, no further $k$ are considered, and that specific group of $k$ hypotheses is identified as containing likely nonnull hypotheses. Viewed as a hypothesis testing procedure for the intersection null hypothesis, Benjamini and Hochberg [4] show that this procedure has level less than or equal to $\alpha$. They also show that procedure controls the FDR, which means, roughly speaking, that, in expectation, a fraction of at least $(1 - \alpha)$ of the rejected null hypotheses should be truly nonnull hypotheses.

How does such a procedure behave in the current setting? We begin by pointing out that Abramovich, Benjamini, Donoho and Johnstone [1] have analyzed the behavior of FDR in exactly the kind of mixture model described above in (1.3)–(1.4), but where $\mu_n$ is calibrated differently from (1.5), and they found an asymptotic minimaxity of FDR in that setting. That is, they considered a situation where one observed data $X_i = \theta_i + Z_i$ and where the $Z_i$ are i.i.d. $N(0,1)$. They supposed that only a fraction $\varepsilon_n$ of the means $\theta_i$ might be nonzero. This is similar to our model above.

As it turns out, for our measure of performance, the behavior of FDR-controlling procedures is not different from that of the maximum $M_n$. To articulate this, consider the detection boundary for the FDR-controlling procedure above. This is the function $\rho_{\mathrm{FDR}}(\beta)$ such that, if $\mu_n = \sqrt{2r\log(n)}$, and if we use a sequence of levels $\alpha = \alpha_n$ tending to 0 slowly enough, then for $r > \rho_{\mathrm{FDR}}(\beta)$ the procedure has power tending to 1 as $n \to \infty$, while for $r < \rho_{\mathrm{FDR}}$ the procedure has power tending to 0.

THEOREM 1.4.

(1.8) $\qquad \rho_{\mathrm{FDR}}(\beta) = \rho_{\mathrm{Max}}(\beta) = (1 - \sqrt{1-\beta})^2, \qquad \tfrac{1}{2} < \beta < 1.$

From our discussion of the behavior of the maximum $M_n$, the FDR-controlling procedure is effective in the range $\beta \in [\tfrac{3}{4}, 1)$, while it is relatively inefficient for $\beta < \tfrac{3}{4}$. Compare Figure 1.

1.5. *Classical methods for combining p-values.* In the literature of meta-analysis [3], one also faces the problem of combining several $p$-values to achieve an overall test of significance. In that literature the component non-null hypotheses are all the same ("homogeneity"), whereas in our discussion they vary widely ("heterogeneity"). A classical approach to combining $p$-values is Fisher's method [13], which refers to $F_n = -2\sum_{1 \le i \le n} \log(p_i)$ as the $\chi^2_n$ distribution. In our setting of extreme heterogeneity, Fisher's method is unable to function well asymptotically:

THEOREM 1.5. *If $\varepsilon_n = n^{-\beta}$, $\beta > \tfrac{1}{2}$ and $\mu_n \le \sqrt{2\log(n)}$, asymptotically $F_n$ is unable to separate $H_1^{(n)}$ and $H_0$.*

1.6. *Relation to goodness-of-fit testing.* Of course, the method we are discussing may be viewed as an application of a goodness-of-fit measure, comparing the empirical distribution of $p$-values to the uniform distribution. As such, it may be compared to many goodness-of-fit procedures where distribution under $H_1$ differs from that of $H_0$ in one tail.

Thus, Anderson and Darling [2] defined a goodness-of-fit measure which involves the maximum of the normalized empirical process. Translated into the current setting, this initially seems very close to HC*. The main difference is that we focus attention near $p = 0$, while Anderson and Darling maximize over $a < p \le b$ with $0 < a < b < 1$. However—important point—all the information needed for discrimination between $H_0^{(n)}$ and $H_1^{(n)}$ is at values of $p$ increasingly close to 0 as $n$ increases. Therefore, statistics based on $a < p < b$ with $a, b$ fixed are dramatically inefficient. Similarly, Borovkov and Sycheva ([7] and [8]) proposed statistics based on the maximum of a normalized empirical process with a general nonlinear normalization $g(p)$ not necessarily $\sqrt{p(1-p)}$. However, the maximum is over $a < p < b$, $0 < a < b < 1$, so our remarks on the Anderson–Darling statistic apply.

Berk and Jones [5] proposed a goodness-of-fit measure which, adapted to the present setting, may be written as

$$\text{BJ}_n^+ = n \cdot \max_{1 \leq i \leq n/2} K^+(i/n, p_{(i)}), \tag{1.9}$$

where $K^+$ is defined by

$$K^+(t, x) = \begin{cases} t \log \dfrac{t}{x} + (1-t) \log \dfrac{1-t}{1-x}, & \text{if } 0 < x < t < 1, \\ 0, & \text{if } 0 \leq t \leq x \leq 1, \\ +\infty, & \text{otherwise,} \end{cases}$$

and is motivated by large deviation theory. Roughly speaking, Lemma 6.4 below shows that $K^+(t,x)$ behaves as $\frac{1}{2}((t-x)^2)/(x(1-x))$, so it should not be surprising that for our measure of performance the detection boundary of $\text{BJ}_n^+$ is the same as that of $\text{HC}^*$. (Note, however, that full justification of the asymptotic claim in [5] has only recently been provided; see [35] for a thorough analysis, which also may shed light on the limiting distribution of $\text{HC}_n^*$.) To articulate our claim about the detection boundary of the Berk–Jones method, define a function $\rho_{\text{BJ}}(\beta)$ such that, if $\mu_n = \sqrt{2r \log(n)}$, and if we use a sequence of levels $\alpha = \alpha_n$ tending to 0 slowly enough, then for $r > \rho_{\text{BJ}}(\beta)$, the procedure has power tending to 1 as $n \to \infty$, while for $r < \rho_{\text{BJ}}$ the procedure has power tending to 0.

THEOREM 1.6.  $\rho_{\text{BJ}}(\beta) \equiv \rho^*(\beta)$, $\frac{1}{2} < \beta < 1$.

However, $\text{HC}^*$ is still better than $\text{BJ}_n^+$ in important ways; we will discuss this in the Appendix, where we prove Theorem 1.6.

1.7. *Generalizations.*  As we have defined it, the higher criticism statistic can obviously be used in a wide variety of situations and there is no need for the $p$-values to be derived from normally distributed $Z$-scores, for example. Consequently, numerous other settings for its deployment can be considered. We have found that, in a wide variety of settings where one has data which are "sparsely nonnull," the $\text{HC}^*$ statistic has an adaptive optimality.

To give the flavor of one of these, we consider a model deriving from the covert communications example mentioned earlier. The problem is one of noncooperative spread-spectrum signal detection. Here one observes $n$ periodogram ordinates $X_i$. In the "covert data absent" case, these represent periodogram ordinates of a white noise, while in the "covert data present" case a small fraction of periodogram ordinates are inflated by the presence of covert signal. The formal model takes the form:

$$H_0 : X_i \stackrel{\text{i.i.d.}}{\sim} \text{Exp}(2), \qquad 1 \leq i \leq n,$$
$$H_1^{(n)} : X_i \stackrel{\text{i.i.d.}}{\sim} (1-\varepsilon)\text{Exp}(2) + \varepsilon \chi_2^2(\delta), \qquad 1 \leq i \leq n.$$

[Here $\text{Exp}(2)$ denotes the exponential distribution with mean 2 and $\chi_2^2(\delta)$ denotes the noncentral chi-squared distribution with noncentrality parameter $\delta$.] Here the data are non-Gaussian, $n$ is large, and the sparsity parameter $\varepsilon = \varepsilon_n = n^{-\beta}$ as before. The strength of the covert signal is measured by the noncentrality parameter $\delta$ of the chi-squared distribution, which we take as $\delta = \delta_n = 2r\log(n)$ for an underlying amplitude parameter $r$ having much the same interpretation as before. (In a cooperative signal detection setting we would know a priori which of the coordinates $X_i$ will exhibit the presence of the covert signal; in the noncooperative case we would not know this.)

We can again apply the principle of higher criticism in this setting, defining $p$-values through the component null hypotheses:

$$p_i = \text{Prob}\{\text{Exp}(2) > X_i\}, \qquad i = 1, \ldots, n.$$

We can also define the detection boundary for this test, as before. This is the function $\rho_{\text{HC,Exp}}(\beta)$ such that, if $\delta_n = 2r\log(n)$, and if we use a sequence of levels $\alpha = \alpha_n$ tending to 0 slowly enough, then, for $r > \rho_{\text{HC,Exp}}(\beta)$, the procedure has power tending to 1 as $n \to \infty$, while, for $r < \rho_{\text{HC,Exp}}$, the procedure has power tending to 0. We can also define the intrinsic detection boundary, the function $\rho^*_{\text{Exp}}(\beta)$ such that if $r < \rho^*_{\text{Exp}}$ the two hypotheses merge asymptotically, while if $r > \rho^*_{\text{Exp}}$ the two hypotheses separate asymptotically.

THEOREM 1.7.

$$(1.10) \quad \rho_{\text{HC,Exp}}(\beta) = \rho^*_{\text{Exp}}(\beta) = \begin{cases} \beta - \frac{1}{2}, & \frac{1}{2} < \beta \leq \frac{3}{4}, \\ (1-\sqrt{1-\beta})^2, & \frac{3}{4} < \beta < 1. \end{cases}$$

In words, the higher criticism statistic achieves the optimal detection region in the $(r, \beta)$ plane; interestingly, this region is the same as we had in the Gaussian case. See Figure 3; also compare to Figure 1.

Other non-Gaussian settings are discussed in Sections 5.1 and 5.2. In each case the higher criticism statistic achieves the (interior of) the optimal detection region.

1.8. *Contents of this paper.* In this paper we establish the key results referred to so far, and then we consider various generalizations. Section 2 develops Theorems 1.1 and 1.2. Section 3 discusses $HC^+$ as a variant of $HC_n^*$ but with better performance in finite samples. Section 4 describes some simulation experiments. Section 5 considers numerous non-Gaussian settings and shows consistently the superiority of higher criticism over Bonferroni and other ideas from the field of multiple comparisons; the proof of Theorem 1.7 is also in this section. The Appendix provides proofs of Theorems 1.3–1.6 as well as Lemma 5.1.

**2. Main results.** Before continuing the narrative flow of the paper, we pause to prove our key results: Theorems 1.1 and 1.2. Later sections will return to the narrative format.

2.1. *Proof of Theorem* 1.1. The idea behind Theorem 1.1 is to simply apply known results from the theory of empirical processes. We recall the normalized uniform empirical process $W_n(t)$ introduced in Section 1.3, and remind the reader that, under $H_0$,

$$HC^* = \max_{0<t<1/2} W_n(t).$$

The normalized empirical process has been studied carefully by a number of authors, and a summary of results can be found in [27], Chapter 16. There, in (16.20), they show that

(2.11) $$\frac{\max_{0<t<\alpha_0} W_n(t)}{\sqrt{2\log\log(n)}} \xrightarrow{p} 1, \qquad n \to \infty,$$

by an argument that depends on the work of Jaeschke [22]. This, in turn, depends on the approximation of $W_n$ by a Brownian bridge and on a result of Darling and Erdös [11], which says that, if $B(t)$ is standard Brownian motion starting at $B(0) = 0$, then

$$\sup_{[1,u]} \frac{B(t)}{\sqrt{t}} \frac{1}{\sqrt{2\log\log u}} \xrightarrow{p} 1, \qquad u \to \infty.$$

Of course, (2.11) implies Theorem 1.1. □

2.2. *Proof of Theorem* 1.2. We begin with a simple observation.

LEMMA 2.1. *Let $X_1, \ldots, X_n$ be i.i.d.* Bernoulli$(\pi_n)$ *and let $a_1, \ldots, a_n$ be a sequence of real numbers. If $n\pi_n \to \infty$ and $a_n/\sqrt{n\pi_n} \to \infty$, then $\lim_{n\to\infty} P(\sum_{i=1}^n (X_i - \pi_n) \leq -a_n) = 0$.*

PROOF. Since $\text{Var}(\sum_{i=1}^{n} X_i) \leq n\pi_n$, use of Chebyshev's inequality yields

$$P\left\{\sum_{i=1}^{n}(X_i - \pi_n) \leq -a_n\right\} \leq \frac{n\pi_n}{a_n^2} \to 0, \qquad n \to \infty. \qquad \square$$

To prove Theorem 1.2, we note that it is enough to show

(2.1) $$\lim_{n\to\infty} P_{H_1^{(n)}}\{\text{HC}_n^* \leq \sqrt{4\log\log(n)}\} = 0.$$

Now, recalling the definition of $V_n(q)$, we have

$$\text{HC}_n^* \geq \sup_{0<q\leq 1} V_n(q).$$

For $0 < q \leq 1$, recall that $p_{n,q} = P\{N(0,1) \geq z_n(q)\}$; we also put

$$p'_{n,q} = P\{(1-\varepsilon_n)N(0,1) + \varepsilon_n N(\mu_n, 1) z_n(q)\}.$$

We now consider two cases. First, suppose that $r > \rho^*(\beta)$ and $r \geq \frac{1}{4}$; then the hypothesis is detectable, if at all, merely by looking at the maximum of the $X_i$. Note that $\sqrt{r} + \sqrt{1-\beta} > 1$. Now, as $V_n(1) \leq \sup_{0<q\leq 1} V_n(q)$ and $V_n(1) =_D (N_n(1) - np_{n,1})/\sqrt{p_{n,1}(1-p_{n,1})}$,

$$P_{H_1^{(n)}}\left\{\sup_{0<q\leq 1} V_n(q) \leq \sqrt{4\log\log(n)}\right\}$$
$$\leq P_{H_1^{(n)}}\{V_n(1) \leq \sqrt{4\log\log(n)}\}$$
$$\leq P_{H_1^{(n)}}\{N_n(1) \leq \sqrt{np_{n,1}}\sqrt{4\log\log(n)} + np_{n,1}\}.$$

Under $H_1^{(n)}$, $N_n(1)$ is a sum of independent Bernoulli($p'_{n,1}$), and by direct calculations

$$p_{n,1} = o(1),$$
$$p'_{n,1} = \frac{1}{\sqrt{4\pi\log(n)}(1-\sqrt{r})} n^{-\beta-(1-\sqrt{r})^2}(1+o(1)).$$

With Lemma 2.1 in mind, let $\pi_n = p'_{n,1}$ and

$$a_n = np'_{n,1} - [\sqrt{np_{n,1}}\sqrt{4\log\log(n)} + np_{n,1}] = np'_{n,1}(1+o(1)),$$

so $n\pi_n \to \infty$ and $a_n/\sqrt{n\pi_n} \to \infty$; the desired result (2.1) in this case follows from Lemma 2.1. In the second case, suppose $r > \rho^*(\beta)$ and $r < \frac{1}{4}$. It will turn out that the hypothesis is detectable based on HC$^*$ but not on the maximum of the $X_i$. Note that $\beta + r < 1$.

$$P_{H_1^{(n)}}\left\{\sup_{0<q\leq 1} V_n(q) \leq \sqrt{4\log\log(n)}\right\}$$
$$\leq P_{H_1^{(n)}}\{V_n(4r) \leq \sqrt{4\log\log(n)}\}$$
$$\leq P_{H_1^{(n)}}\{N_n(4r) \leq \sqrt{np_{n,4r}}\sqrt{4\log\log(n)} + np_{n,4r}\}.$$

Under $H_1^{(n)}$ $N_n(4r)$ is a sum of $n$ independent Bernoulli$(p'_{n,4r})$, and

$$p_{n,4r} = \frac{1}{4\sqrt{\pi r \log(n)}} n^{-4r}(1+o(1)),$$

$$p'_{n,4r} = p_{n,4r} + \frac{1}{\sqrt{4\pi r \log(n)}} n^{-(\beta+r)}(1+o(1)).$$

With Lemma 2.1 in mind, let $\pi_n = p'_{n,4r}$ and

$$a_n = n[p'_{n,4r} - p_{n,4r}] - \sqrt{np_{n,4r}}\sqrt{4\log\log(n)}.$$

Direct calculation shows

$$\frac{a_n}{\sqrt{n\pi_n}} = \begin{cases} O(n^{r-(\beta-1/2)}/(\log(n))^{1/4}), & \text{if } \beta > 3r, \\ O(n^{1/2[1-(\beta+r)]}/(\log(n))^{1/4}), & \text{if } \beta \leq 3r, \end{cases}$$

so $n\pi_n \to \infty$, $a_n/\sqrt{n\pi_n} \to \infty$; the desired result (2.1) follows for this case from Lemma 2.1. □

**3. A refinement.** In the most challenging cases, where $0 < r < \frac{1}{4}$, the analysis above tells us that the real information about the presence of $H_1^{(n)}$ is located away from the extreme values, and so, perhaps unexpectedly, the few smallest $p_i$ are not relevant to the detection problem.

This suggests that we will still be able to reach the full interior of the detection region if we work with a modified statistic involving the maximum over all $p$-values greater than or equal to $1/n$, which is

$$\text{HC}_n^+ = \max_{1 < i \leq n/2,\ p_{(i)} \geq 1/n} \sqrt{n}\,[i/n - p_{(i)}]/\sqrt{p_{(i)}(1-p_{(i)})}.$$

Of course, this restriction is not necessary from the viewpoint of asymptotic analysis, since our results show that $\text{HC}_n^*$ is effective without any adjustment. However, our experience with moderate-sized samples shows this adjustment to be quite valuable.

The empirical process viewpoint is helpful. Under $H_0$ $\text{HC}_n^*$ behaves as

(3.1) $$\sup_{0<t<1/2} W_n(t).$$

For comparison, $\mathrm{HC}_n^+$ behaves as

$$\sup_{1/n < t < 1/2} W_n(t). \tag{3.2}$$

The limitation of the range in (3.2) is important. The quantity in (3.1) can be very much dominated by behavior in the vicinity of 0, or, more particularly, by the smallest observation. Since the smallest observation (smallest $p$-value) is not where the information for detection resides, this emphasis is misplaced.

Speaking quantitatively, the statistic can have a heavy-tailed distribution under the null hypothesis, as we will show in a moment. Since, ordinarily, one does not want to use test statistics with heavy-tailed distributions— because tests at stringent levels would have low power in small samples—the quantitative motivation for $\mathrm{HC}_n^+$ is clear.

The heavy tails under the null hypothesis come from the following effect. Recall that, under $H_0$, the minimum $p$-value has a distribution that is $\mathrm{Exp}(1/n)$. It follows that the first component in the maximum defining $\mathrm{HC}_n^*$,

$$\mathrm{HC}_{n,1}^* = \sqrt{n}\,[1/n - p_{(1)}]/\sqrt{p_{(1)}(1-p_{(1)})},$$

has the asymptotic distribution

$$\mathrm{HC}_{n,1}^* \stackrel{D}{\Longrightarrow} \frac{1}{\sqrt{E}} - \sqrt{E},$$

where $E$ is exponentially distributed with mean 1. Now, for large $t$,

$$P\left(\frac{1}{\sqrt{E}} - \sqrt{E} > t\right) = P(\sqrt{E} < (\sqrt{t^2+4}-t)/2) \sim \frac{1}{t^2}.$$

Hence, $\mathrm{HC}_{n,1}^*$ has "heavy tails" under $H_0$.

Numerical simulations show that unusually large values of $\mathrm{HC}_n^*$ under the null hypothesis are most frequently caused by $\mathrm{HC}_{n,1}^*$, suggesting that if we restrict the range of the maximum as in $\mathrm{HC}_n^+$ this "heavy tail" is thinned out considerably. Numerical experiments support this analysis. Adapting material from [12], or from [27], page 600, it seems we have the following asymptotic law for $\mathrm{HC}_n^+$:

$$b_n \mathrm{HC}_n^+ - c_n \to E_v^2,$$

where

$$b_n = \sqrt{2\log\log(n)}, \qquad c_n = 2\log\log(n) + \tfrac{1}{2}\log(\log\log(n)) - \tfrac{1}{2}\log(4\pi),$$

and the c.d.f. of $E_v^2$ is $\exp(-2\exp(-x))$. Experiments show that this is a fairly accurate approximation for moderate $n$. The same form of limiting distribution holds for $\mathrm{HC}_n^*$, of course, but it is empirically a poor fit.

**4. Simulation.** We have conducted a small-scale empirical study of the performance of $HC_n^*$ and $HC_n^+$. Our idea was to select a few interesting $(r,\beta)$ pairs in the detectable region—above, but close to the boundary, to create samples from null and alternative—and to study the behavior of the proposed statistics.

We took $n = 10^6, \varepsilon = \frac{1}{1000}$ and $\mu = \sqrt{2 \cdot 0.15 \cdot \log(n)} \approx 2.04$ for the following experiment:

1. Draw $n$ samples from $N(0,1)$ to represent $H_0$; calculate $HC^*$ and $HC^+$.
2. Replace 1000 of the previous samples by the same number of samples from $N(\mu,1)$ to represent $H_1^{(n)}$; calculate $HC^*$ and $HC^+$.
3. Repeat 1 and 2 100 times and make histograms of the simulated $HC^*$ and $HC^+$.

See Figure 2 for the results. As can be seen from the theory above, detectability requires increasingly large samples as one approaches the detection boundary. In fact, depending on how close one goes to the boundary, the required sample sizes can become enormous. To have an idea of how large $n$ may have to be, we consider the case $r < \frac{1}{4}$, $\frac{1}{2} < \beta < \frac{3}{4}$. We recall that

$$(4.1) \quad EV_n(4r) \begin{cases} = 0, & \text{under } H_0, \\ \sim n^{r-(\beta-1/2)}/\sqrt[4]{\pi r \log(n)}, & \text{under } H_1^{(n)}. \end{cases}$$

Numerically, we calculate the values of $EV_n(4r)$ under $H_1^{(n)}$ for various $n$ with $r = 0.1$, $\beta = \frac{1}{2}$ and $r = 0.05$, $\beta = \frac{1}{2}$, respectively; the results are summarized in Table 1, together with the values of $\sqrt{2\log\log(n)}$ for comparison; notice the variance of $V_n(4r)$ is roughly 1. This, of course, raises the issue of whether our theory adequately describes practice in small samples and, in particular, whether the motivating examples of the Introduction really offer credible scenarios for deployment of these results. We leave such discussion to another occasion.

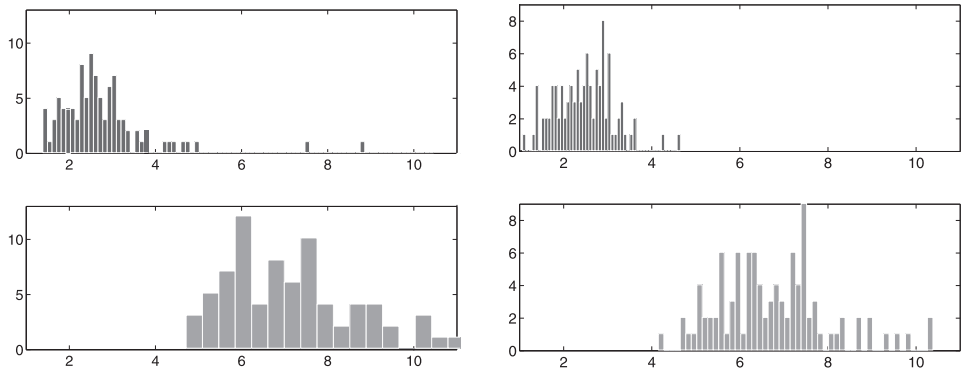

FIG. 2. *Histograms for* $HC^*$ *and* $HC^+$. *Top row: Behavior under* $H_0$. *Bottom row: Behavior under* $H_1^{(n)}$. *Left column:* $HC^*$. *Right column:* $HC^+$. *Here* $n = 10^6$, $\mu = \sqrt{0.3\log(n)} \approx 2.04$, $\varepsilon = 10^{-3}$.

TABLE 1
$EV_n(4r)$ under $H_1^{(n)}$ for various $n$. The parameters $(r, \beta)$ are
$r = 0.1$, $\beta = \frac{1}{2}$ and $r = 0.05$, $\beta = \frac{1}{2}$, respectively; values of $\sqrt{2 \log \log(n)}$
are also included for comparison

| $n$ | $10^6$ | $10^7$ | $10^8$ | $10^9$ | $10^{10}$ |
|---|---|---|---|---|---|
| $\sqrt{2 \log \log(n)}$ | 2.2916 | 2.3579 | 2.4139 | 2.4622 | 2.5046 |
| $EV_n(4r)$ $r = 0.1$, $\beta = \frac{1}{2}$ | 2.7582 | 3.3411 | 4.0680 | 4.9728 | 6.0976 |
| $r = 0.05$, $\beta = \frac{1}{2}$ | 1.6439 | 1.7748 | 1.9259 | 2.0982 | 2.2931 |

The large-sample issue also raises the question of how to efficiently sample from the normal distribution with very large $n$. Our approach goes as follows.

1. Pick a small number $\varepsilon$, such as $\varepsilon = 10^{-3}$ or $\varepsilon = 10^{-6}$.
2. Simulate samples from uniform distribution at quantiles greater than $1 - \varepsilon$.

   - Sample a number $K$ from the Poisson distribution with mean $n\varepsilon$.
   - Generate $K$ samples $(U_1, \ldots, U_K)$ from the uniform distribution on $(1 - \varepsilon, 1)$.
   - Generate $K$ samples of $(z_1, \ldots, z_K)$ by letting

$$(4.2) \quad y_i = -2 \log(1 - U_i) - \log 2\pi, \qquad z_i = \sqrt{y_i - \log y_i}, \qquad 1 \leq i \leq K.$$

Approximately, the $(z_1, \ldots, z_K)$ can be viewed as if a sample of size $n$ had been taken from the Normal distribution, and then only the $(1 - \varepsilon)n$ largest sample values were retained. Compared to brute force, the algorithm requires only $\varepsilon n$ flops, rather than $n$ flops. Obviously, the accuracy of the approximation in (4.2) depends on how small $\varepsilon$ is; the smaller the $\varepsilon$, the more accurate the approximation, and thus the more accurate the simulation.

**5. Other settings.** The principle of higher criticism can be applied in a much wider series of situations than the Gaussian model studied so far, as we now show.

5.1. *Chi-Squared.* Let $\chi_\nu^2(\delta)$ denote the usual chi-squared distribution with $\nu$ degrees of freedom and noncentrality parameter $\delta$. Consider the problem of testing between these hypotheses:

$$(5.1) \qquad H_0 : X_i \stackrel{\text{i.i.d.}}{\sim} \chi_\nu^2(0), \qquad\qquad 1 \leq i \leq n,$$

$$(5.2) \qquad H_1^{(n)} : X_i \stackrel{\text{i.i.d.}}{\sim} (1 - \varepsilon)\chi_\nu^2(0) + \varepsilon\chi_\nu^2(\delta), \qquad 1 \leq i \leq n.$$

Here $\varepsilon = \varepsilon_n = n^{-\beta}$ as before. Owing to the representation of chi-squared r.v.'s as sums of squares of standard normals, the problem can be rewritten in terms of arrays of normals $z_{ij} \stackrel{\text{i.i.d.}}{\sim} N(0,1)$:

$$(5.3) \quad H_0 : X_i \stackrel{\text{i.i.d.}}{\sim} \sum_{j=1}^{\nu} z_{ij}^2, \qquad 1 \le i \le n,$$

$$(5.4) \quad H_1^{(n)} : X_i \stackrel{\text{i.i.d.}}{\sim} (1-\varepsilon) \sum_{j=1}^{\nu} z_{ij}^2 + \varepsilon \left[ (z_{i1} + \delta)^2 + \sum_{j=2}^{\nu} z_{ij}^2 \right], \qquad 1 \le i \le n,$$

where $z_{ij} \stackrel{\text{i.i.d.}}{\sim} N(0,1)$. Roughly speaking, a small fraction of the normals may have nonzero mean, and we must base our decisions on sums of squares rather than on the normals themselves. The problem is obviously related to the Gaussian problem discussed above; only the observations with nonzero means occur in groups, and while the nonzero means are not equal within a group, the sums of squares of those means within a group must be equal.

Now, obviously, if $\nu = 1$, then we are equivalently dealing with squares of individual normals, and so are just seeing a two-sided variant of the original one-sided normal testing problem considered so far. As it turns out, all the detection boundary and attainability results for the two-sided normal problem are the same as in the one-sided case.

If $\nu = 2$, we view this as modeling the covert communications problem of Section 1.5. Indeed, the real and imaginary components of the discrete Fourier transform of Gaussian white noise are normal and independent; the sum of squares of those two components is just the periodogram. In a frequency where there is no signal, only Gaussian noise, the periodogram has a $\chi_2^2(0)$ [or $\text{Exp}(2)$] distribution, which is precisely the exponential mentioned earlier; while, in the signal-present case, the periodogram has a $\chi_2^2(\delta_k)$ distribution, where $\delta_k$ is the signal energy at that frequency.

If $\nu > 2$, we can think of agricultural trials, where a treatment is attempted with replications and in many different blocks. The alternative hypothesis is that, in a few special blocks, the treatment has a substantial effect, but in most blocks it has no effect.

As it turns out, for any fixed, constant number of degrees of freedom, the results will be similar. Let the noncentrality parameter obey

$$\delta = \delta_n = 2r \log(n), \qquad 0 < r < 1.$$

In this problem, we again have an $(r, \beta)$ plane, in which there is a region of detectability and a detection boundary. There are two key facts:

- The optimal detection boundary $\rho^*_{\chi_2^2}(\beta)$ is the same as in the Gaussian case:

$$\rho^*_{\chi_2^2}(\beta) = \begin{cases} (1 - \sqrt{1-\beta})^2, & \frac{3}{4} \le \beta \le 1, \\ \beta - \frac{1}{2}, & \frac{1}{2} < \beta < \frac{3}{4}. \end{cases}$$

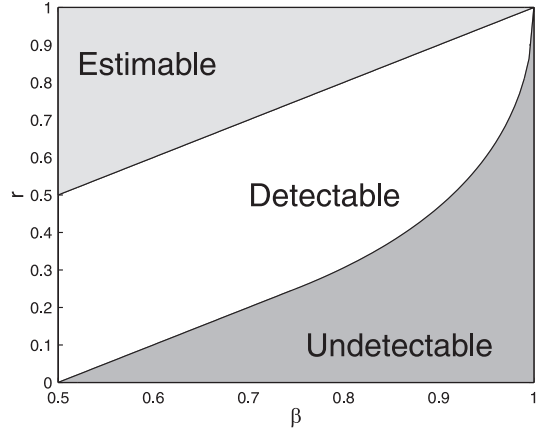

FIG. 3. *Three regions in the $r$–$\beta$ plane for the $\chi^2$ model (5.1)–(5.2). The detection boundary separates the detectable region from the undetectable region. For the estimable region, it is possible not only to detect the presence of nonzero means, but also to estimate those means.*

This is proven in [23].
- $HC_n^*$ is able to separate the null and alternative hypotheses throughout the interior of the detection region, and thus its detection boundary $\rho_{HC,\chi_2^2}$ obeys

$$\rho_{HC,\chi_2^2}(\beta) = \rho^*_{\chi_2^2}(\beta), \qquad \tfrac{1}{2} < \beta \leq 1.$$

The key ideas are illustrated in Figure 3.

The analysis supporting the performance of $HC_n^*$ is as follows. Define $x_n^2(q) = 2q\log(n)$. We then have

(5.5)
$$P\{\chi_\nu^2(0) > x_n^2(q)\} \sim L_n \cdot n^{-q},$$
$$P\{\chi_\nu^2(\delta_n) > x_n^2(q)\} \sim L_n \cdot n^{-(\sqrt{q}-\sqrt{r})^2}, \qquad r < q \leq 1.$$

The left-hand relation is proved as follows:

$$P\{\chi_\nu^2(0) > x_n^2(q)\} = \frac{2^{1-\nu/2}}{\Gamma(\frac{\nu}{2})} \int_{\sqrt{2q\log(n)}}^{\infty} \rho^{\nu-1} e^{-\rho^2/2} d\rho$$
$$= \frac{(q\log(n))^{\nu/2-1}}{\Gamma(\frac{\nu}{2})} n^{-q}(1+o(1)).$$

It takes a little more effort to check that the second equality is also true; this follows from the following lemma.

LEMMA 5.1. *Let $0 < r < q < 1$ and $\delta_n = 2q\log(n)$. Then, as $n \to \infty$,*
$$P\{\chi_\nu^2(\delta_n) \geq x_n^2(q)\}$$

$$= \frac{1}{\sqrt{2\pi \log(n)}} \left(\frac{r}{q}\right)^{(1-\nu)/4} \frac{1}{\sqrt{2q} - \sqrt{2r}} n^{-(\sqrt{q}-\sqrt{r})^2}(1+o(1)).$$

The proof of Lemma 5.1 is given in the Appendix. Once the relations (5.5) are available, the analysis proceeds exactly as in the earlier Gaussian case. The most informative values of $q$ (resp. $x^2$) are at

$$\text{if } r < \tfrac{1}{4}, \quad \text{then } q = 4r \text{ or } x^2 \approx 4\delta_n,$$
$$\text{if } r \geq \tfrac{1}{4}, \quad \text{then } q = 1 \text{ or } x^2 \approx 2\log(n),$$

correspondingly.

5.2. *Generalized Gaussian (Subbotin) distribution.* The generalized Gaussian (Subbotin) distribution $\text{GN}_\gamma(\mu)$ has density function

$$\frac{1}{C_\gamma} \exp\left(-\frac{|x-\mu|^\gamma}{\gamma}\right), \qquad C_\gamma = 2\Gamma\left(\frac{1}{\gamma}\right)\gamma^{1/\gamma - 1}.$$

This class of densities was introduced by Subbotin [29]; see [24], page 195. It has many uses in Bayesian analysis; see [9], page 157, who cite several earlier references. The Gaussian is, of course, the special case $\gamma = 2$. The case $\gamma = 1$ corresponds to the double exponential (Laplace) distribution, a well-understood and widely used distribution. The case $\gamma < 1$ is of interest in image analysis of natural scenes, where it has been found that wavelet coefficients at a single scale can be modeled as following a Subbotin distribution with $\gamma \approx 0.7$ [28]. This suggests that various problems of image detection, such as in watermarking and steganography, could reasonably use the model above.

A natural generalization of (1.3)–(1.4) is the following:

(5.6) $\qquad H_0 : X_i \stackrel{\text{i.i.d.}}{\sim} \text{GN}_\gamma(0), \qquad\qquad 1 \leq i \leq n,$

(5.7) $\qquad H_1^{(n)} : X_i \stackrel{\text{i.i.d.}}{\sim} (1-\varepsilon)\text{GN}_\gamma(0) + \varepsilon \text{GN}_\gamma(\mu), \qquad 1 \leq i \leq n.$

Here we choose the calibrations

$$\varepsilon_n = n^{-\beta}, \qquad \mu = \mu_{\gamma,n} = (\gamma r \log(n))^{1/\gamma}, \qquad \tfrac{1}{2} < \beta < 1, \ 0 < r < 1.$$

First we will discuss the case $\gamma > 1$. In this range the number of nonzero means is too small to be noticeable in any sum which is of expectation of order $n$; if $r$ is not large, we cannot expect a visible effect on the upper extreme. In short, this detection problem will be a difficult problem.

For $\gamma > 1$, we showed in [23] that the detection boundary is defined as

(5.8) $\rho^*_\gamma(\beta) = \begin{cases} (2^{1/(\gamma-1)} - 1)^{\gamma-1}(\beta - \tfrac{1}{2}), & \tfrac{1}{2} < \beta \leq 1 - 2^{-\gamma/(\gamma-1)}, \\ (1 - (1-\beta)^{1/\gamma})^\gamma, & 1 - 2^{-\gamma/(\gamma-1)} \leq \beta < 1. \end{cases}$

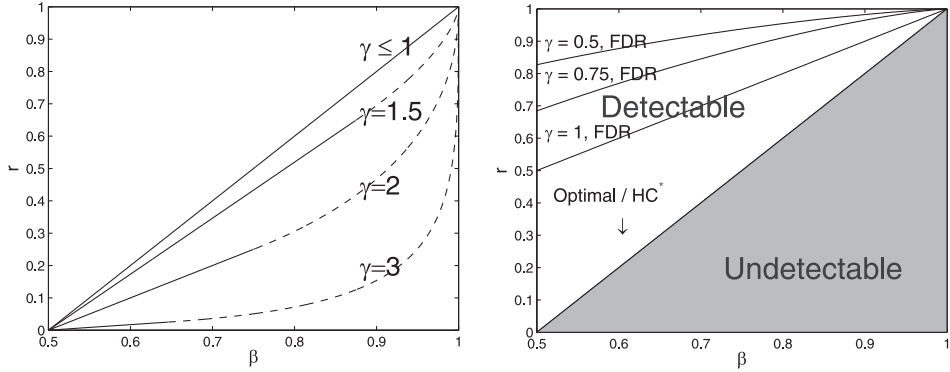

FIG. 4. Left: *Detection boundaries of the $r$–$\beta$ plane for model* (5.6)–(5.7) *for $\gamma \leq 1$ and $\gamma = 1.5$, 2, 3 from top to bottom. Solid parts of the curves are line segments.* Right: *The common detection boundary for all $0 < \gamma \leq 1$ which separates the detectable region from the undetectable region. Three curves from top to bottom correspond to the detection boundaries of the Bonferroni method with $\gamma = \frac{1}{2}$, $\gamma = \frac{3}{4}$ and $\gamma = 1$.*

Now we discuss the case $\gamma \leq 1$. Jin [23] showed that there is a threshold effect for the likelihood ratio test, but the detection boundary is quite different, and, surprisingly, it can be described in terms of $(r, \beta)$ independent of $\gamma$:

(5.9) $$\rho_\gamma^*(\beta) = 2(\beta - \tfrac{1}{2}), \qquad \tfrac{1}{2} < \beta < 1.$$

We have the following result.

THEOREM 5.1. *Consider applying higher criticism to the p-values $p_i = P\{GN_\gamma(0) > X_i\}$, $i = 1, \ldots, n$, in the setting just described. Then the detection boundary $\rho_{HC,\gamma}$ for this procedure is the same as the efficient detection boundary:*

$$\rho_{HC,\gamma}(\beta) = \rho_\gamma^*(\beta), \qquad \tfrac{1}{2} < \beta < 1.$$

The basic phenomena are depicted in Figure 4.

The analysis can be made very similar to the normal case. Specifically, introduce the notation $z_{\gamma,n}(q) = (\gamma q \log(n))^{1/\gamma}$ and let $\pi_{\gamma,n,q} = P\{GN_\gamma(0) > z_{\gamma,n}(q)\}$, $0 < q < 1$. Note that, when $\gamma = 2$, $z_{\gamma,n}(q) \equiv z_n(q)$ and $\pi_{\gamma,n,q} \equiv p_{n,q}$. We have

$$\max_{0 < t < 1/2} W_n(t) = \max_{0 < q < \infty} W_n(\pi_{\gamma,n,q}).$$

Similarly, let $N_{\gamma,n}(q)$ count the observations exceeding $z_{\gamma,n}(q)$:

$$N_{\gamma,n}(q) = \#\{i : X_i \geq z_{\gamma,n}(q)\}$$

and also
$$V_{\gamma,n}(q) = \frac{N_{\gamma,n}(q) - n\pi_{\gamma,n,q}}{\sqrt{n\pi_{\gamma,n,q}(1 - \pi_{\gamma,n,q})}}.$$

By arguments which are obvious at this point,
$$P\{\text{GN}_\gamma(\mu_{\gamma,n}) > z_{\gamma,n}(q)\} = L_n \cdot n^{-(q^{1/\gamma} - r^{1/\gamma})^\gamma},$$
$$P\{\text{GN}_\gamma(0) > z_{\gamma,n}(q)\} = L_n \cdot n^{-q}, \qquad r < q \leq 1.$$

It follows that, under the alternative $H_1^{(n)}$, we have
$$EV_{\gamma,n}(q) = L_n \cdot \frac{n^{1-\beta} \cdot n^{-(q^{1/\gamma} - r^{1/\gamma})^\gamma}}{\sqrt{n \cdot n^{-q}}}$$
$$= L_n \cdot n^{[(1+q)/2 - \beta - (q^{1/\gamma} - r^{1/\gamma})^\gamma]},$$

while under the null $EV_{\gamma,n}(q) = 0$. The most informative value of $q$ will optimize the growth rate of $n^{[(1+q)/2 - \beta - (q^{1/\gamma} - r^{1/\gamma})^\gamma]}$ to $\infty$.

For the case $\gamma > 1$, define
$$r_\gamma = (1 - 2^{-1/(\gamma-1)})^\gamma.$$

The most informative value of $q$ satisfies if $r < r_\gamma$, then $q = \dfrac{r}{r_\gamma}$ and $EV_n(q) = L_n \cdot n^{[(2^{1/(\gamma-1)} - 1)^{1-\gamma} \cdot r - (\beta - 1/2)]}$, if $r \geq r_\gamma$, then $q = 1$ and $EV_n(q) = L_n \cdot n^{[(1-\beta) - (1 - r^{1/\gamma})^\gamma]}$.

For the case $0 < \gamma \leq 1$, the story is quite different, and the main reason is that
$$\frac{1+q}{2} - \beta - (q^{1/\gamma} - r^{1/\gamma})^\gamma,$$
as a function of $q$, is strictly decreasing for any fixed $0 < \gamma \leq 1$, so the most informative place to look is at

$$q = r \quad \text{or, equivalently, at } x \approx \mu.$$

Notice that, under $H_0$, HC$^*$ behaves the same as in the normal case. Under $H_1^{(n)}$ the above analysis shows the behavior at the most informative place. We can argue exactly as in the proof of Theorem 1.2. The growth of $EV_{\gamma,n}(q)$ easily surpasses the $\sqrt{4 \log \log(n)}$ threshold, and the result follows.

There are some interesting points here.

First, the detection boundary for all the cases where $\gamma \leq 1$ looks like the limit of the boundaries for $\gamma > 1$ as $\gamma \to 1$. Second, the most informative place to look, for the case $\gamma > 1$, is at
$$x \approx \frac{1}{1 - 2^{-1/(\gamma-1)}} \mu,$$

the coefficient $1/(1 - 2^{-1/(\gamma-1)}) \to 1$ as $\gamma \to 1$; in comparison, for the case $\gamma \leq 1$, the most informative place to look is at $x = \mu$.

Third, it is interesting to notice that, for the case $\gamma \leq 1$, the best that either the maximum or the FDR-controlling methods can obtain is

(5.10) $$r > (1 - (1-\beta)^{1/\gamma})^\gamma.$$

This is strictly above the detection boundary as defined in (5.9) for any $\frac{1}{2} < \beta < 1$, while, in comparison, higher criticism can obtain the full interior of the region of detectability, for all $\gamma$. Fourth, notice that the performance of the maximum or FDR-controlling methods worsens compared to $HC^*$ when $\gamma \to 0$. The best that the maximum or FDR-controlling methods can do when $\gamma \approx 0$ is to detect for $r > 1$, while higher criticism is able to detect for $r > 2\beta - 1$, $\frac{1}{2} < \beta < 1$, independent of $\gamma$; the superiority of $HC^*$ can be seen most prominently for the case $\beta \approx \frac{1}{2}$, $\gamma \approx 0$, in which $HC^*$ is able to detect for $r > 2\beta - 1 \approx 0$, while the maximum or FDR-controlling methods are able to detect only for $r > 1$. Compare Figure 4.

## APPENDIX: PROOFS

### A.1. Proof of Theorem 1.3.

LEMMA A.2. *If $z_n \sim \text{Binomial}(n, \pi_n)$ and $\pi_n \to 0$, $n\pi_n \to \infty$, then $P\{z_n \geq 1\} \to 1$.*

PROOF.
$$P\{z_n = 0\} = (1 - \pi_n)^n = e^{-n\log(1-\pi_n)} \to 0. \qquad \square$$

PROOF OF THEOREM 1.3. When $r > \rho^+(\beta)$ or, equivalently, $1 - \beta > (1 - \sqrt{r})^2$, we can pick a constant $c > 0$ depending only on $(r, \beta)$ such that
$$1 - \beta > (\sqrt{1+c} - \sqrt{r})^2.$$

To prove Theorem 1.3, it is sufficient to prove

(A.11) $\quad P_{H_1^{(n)}}\{M_n \geq \sqrt{2(1+c)\log(n)}\} \to 1 \quad$ as $n \to \infty$.

Let
$$N(c) = \#\{i : X_i \geq \sqrt{2(1+c)\log(n)}\}.$$

Then, under $H_1^{(n)}$ $N(c) \sim \text{Binomial}(n, q_{n,c})$, where
$$q_{n,c} = P\{(1-\varepsilon_n)N(0,1) + \varepsilon_n N(\mu_n, 1) \geq \sqrt{2(1+c)\log(n)}\}$$
$$= L_n \cdot n^{-\beta - (\sqrt{1+c} - \sqrt{r})^2}.$$

Notice that, as $n \to \infty$, $q_{n,c} \to 0$ and $nq_{n,c} \to \infty$. So letting $\pi_n = q_{n,c}$ in Lemma A.1, we have

$$P_{H_1^{(n)}}\{N(c) \geq 1\} \to 1 \qquad \text{as } n \to \infty.$$

Finally, the desired result (A.11) follows from

$$P_{H_1^{(n)}}\{M_n \geq \sqrt{2(1+c)\log(n)}\,\} = P_{H_1^{(n)}}\{N(c) \geq 1\}. \qquad \square$$

**A.2. Proof of Theorem 1.4.**

LEMMA A.3. *For constants* $\frac{1}{2} < \delta < 1$, $0 < a < 1$,

$$\sup_{\{t\,:\,t > n^{-\delta}\}} P\{\text{Binomial}(n,t) \geq n \cdot t/a\} \leq 2e^{-c(a)n^{1-\delta}}, \qquad n \to \infty,$$

*where* $c(a) = 1 - (\log a + 1)/a$.

PROOF. By noticing that

$$\sup_{\{t > n^{-\delta}\}} P\{\text{Binomial}(n,t) \geq n \cdot t/a\} \leq P\bigg\{\sup_{t > n^{-\delta}} \frac{\text{Binomial}(n,t)}{nt} \geq \frac{1}{a}\bigg\},$$

Lemma A.3 follows directly from [34], Lemma 1. $\square$

Now for any $\frac{1}{2} < \delta < 1$ introduce statistics:

$$F_1^\delta = \min_{\{i\,:\,p_{(i)} > n^{-\delta}\}} \frac{p_{(i)}}{i/n}, \qquad F_2^\delta = \min_{\{i\,:\,p_{(i)} \leq n^{-\delta}\}} \frac{p_{(i)}}{i/n}.$$

LEMMA A.4. *For any constants* $\frac{1}{2} < \delta < 1$, $0 < a < 1$, *if* $r < \rho_{\text{FDR}}(\beta)$, *then, as* $n \to \infty$,

(A.12) $\qquad P_{H_0}\{F_1^\delta \leq a\} \to 0, \qquad P_{H_1^{(n)}}\{F_1^\delta \leq a\} \to 0.$

PROOF. Recall that $nF_n(t) = \sum_i \mathbb{1}_{\{p_i \leq t\}}$, $0 < t < 1$. We have $nF_n(t) \sim \text{Binomial}(n,t)$ under $H_0$ and $nF_n(t) \sim \text{Binomial}(n, \pi(n,t))$ under $H_1^{(n)}$, where $\pi(n,t) = P_{H_1^{(n)}}\{p_i \leq t\} \geq t$; since $r < \rho_{\text{FDR}}(\beta)$ we also have

$$a_n \stackrel{\triangle}{=} \sup_{\{t > n^{-\delta}\}} \pi(n,t)/t \to 1.$$

Observe that $i/n = F_n(p_{(i)})$, so

$$\frac{p_{(i)}}{i/n} \leq a \quad \Longleftrightarrow \quad F_n(p_{(i)}) \geq p_{(i)}/a;$$

by Lemma A.3, the desired result in (A.12) follows from

$$P_{H_0}\{F_1^\delta \leq a\} \leq n \cdot \sup_{\{t: t > n^{-\delta}\}} P\{\text{Binomial}(n, t) \geq n \cdot t/a\} \to 0,$$

$$P_{H_1^{(n)}}\{F_1^\delta \leq a\} \leq n \cdot \sup_{\{t: t > n^{-\delta}\}} P\{\text{Binomial}(n, \pi(n, t)) \geq n \cdot t/a\} \to 0. \quad \square$$

PROOF OF THEOREM 1.4. For the FDR-controlling procedure,

$$\text{Reject if and only if } \min_{1 \leq i \leq n} \frac{p_{(i)}}{i/n} \leq h(n, \alpha_n),$$

where $h(n, \alpha_n) \leq a < 1$ is any given critical value. Since the attainability of the FDR-controlling procedure is as good as the maximum or Bonferroni method, all we need to prove is that, for $(r, \beta)$ in the region $\rho^*(\beta) < r < (1 - \sqrt{1-\beta})^2$, the FDR-controlling method totally fails or

(A.13) $\quad P_{H_0}\{\text{Reject } H_0\} + P_{H_1^{(n)}}\{\text{Accept } H_0\} \to 1 \quad \text{as } n \to \infty.$

Now, under $H_1^{(n)}$ we break $\{1, 2, \ldots, n\}$ into two sets $A_1^{(n)}$ and $A_2^{(n)}$, where

$$i \in A_1^{(n)} \quad \text{if } X_i \text{ is sampled from } N(0, 1),$$

$$i \in A_2^{(n)} \quad \text{if } X_i \text{ is sampled from } N(\mu_n, 1).$$

Introduce an event:

$$E_n^{\delta_0} = \{p_{(i)} \leq n^{-\delta_0} \text{ for some } i \in A_2^{(n)}\}.$$

Since $r < (1 - \sqrt{1-\beta})^2$, or $\sqrt{r} + \sqrt{1-\beta} < 1$, we can choose $\delta_0$ to be close enough to 1 such that $P_{H_1^{(n)}}(E_n^{\delta_0}) \to 0$. Notice that

$$\{F_2^{\delta_0} | H_0\} \stackrel{D}{=} \{F_2^{\delta_0} | (H_1^{(n)}, (E_n^{\delta_0})^c)\}$$

and

$$P_{H_1^{(n)}}\{F_2^{\delta_0} < h(n, \alpha_n)\}$$
$$= (1 - P_{H_1^{(n)}}(E_n^{\delta_0})) P_{H_1^{(n)}}\{F_2^{\delta_0} < h(n, \alpha_n) | (E_n^{\delta_0})^c\}$$
$$+ P_{H_1^{(n)}}(E_n^{\delta_0}) P_{H_1^{(n)}}\{F_2^{\delta_0} < h(n, \alpha_n) | E_n^{\delta_0}\},$$

so $P_{H_0}\{F_2^{\delta_0} < h(n, \alpha_n)\} - P_{H_1^{(n)}}\{F_2^{\delta_0} < h(n, \alpha_n)\} \to 0$.

Finally, by Lemma A.4

$$|P_{H_0}\{\text{Reject}\} - P_{H_1^{(n)}}\{\text{Reject}\}|$$
$$\leq |P_{H_0}\{F_2^{\delta_0} < h(n, \alpha_n)\} - P_{H_1^{(n)}}\{F_2^{\delta_0} < h(n, \alpha_n)\}|$$

$$+ |P_{H_0}\{F_1^{\delta_0} < h(n,\alpha_n)\}| + |P_{H_1^{(n)}}\{F_1^{\delta_0} < h(n,\alpha_n)\}|$$
$$\to 0,$$

and the desired result in (A.13) follows. □

**A.3. Proof of Theorem 1.5.** Under $H_0$, $F_n \sim \chi^2_{2n}$, and $EF_n = 2n$, $\text{Var}(F_n) = 4n$. Under $H_1^{(n)}$ with $\varepsilon_n = n^{-\beta}$, $\mu_n = \sqrt{2r\log(n)}$, $\frac{1}{2} < \beta < 1$, $0 < r \leq 1$, direct calculations show that $2n \leq EF_n = 2n[1 + O(\varepsilon_n L_n)]$, $\text{Var}(F_n) = 4n[1 + O(\varepsilon_n L_n)]$; since $\beta > \frac{1}{2}$, the conclusion follows.

**A.4. Proof of Theorem 1.6.** To prove Theorem 1.6 we need the following lemma.

LEMMA A.5. (i) *For* $0 < x < t \leq \frac{1}{2}$,

(A.14) $$K^+(t,x) \leq \frac{1}{2}\frac{(t-x)^2}{x(1-x)}.$$

(ii) *Let* $x = x(t)$ *obey* $0 < x < t < 1$. *We have, as* $t \to 0$,

(A.15) $$K^+(t,x) = \begin{cases} \dfrac{1}{2}\dfrac{(t-x)^2}{x(1-x)}\left(1 + O\left(t + \dfrac{t}{x} - 1\right)\right), & \text{if } \dfrac{t}{x} \to 1, \\ t\log\dfrac{t}{x}(1+o(1)), & \text{if } \dfrac{t}{x} \to \infty. \end{cases}$$

PROOF. (i) Letting $t = sx$, it is sufficient to prove that, for fixed $0 < x < \frac{1}{2}$ and for $1 \leq s \leq \frac{1}{2x}$,

(A.16) $$s\log s + \left(\frac{1}{x} - s\right)\log\left(\frac{1-sx}{1-x}\right) \leq \frac{1}{2}\frac{(s-1)^2}{1-x}.$$

To prove (A.16), set

$$f(s) = s\log s + \left(\frac{1}{x} - s\right)\log\left(\frac{1-sx}{1-x}\right) - \frac{1}{2}\frac{(s-1)^2}{1-x};$$

direct calculations show that $f(1) = 0$, $f'(1) = 0$ and

$$f''(s) = \frac{(1-s)(1-(s+1)x)}{s(1-x)(1-xs)}, \qquad 1 \leq s \leq \frac{1}{2x}.$$

Notice that when $s \geq 1$, $1 - (s+1)x \geq 1 - 2sx$, so, for any fixed $0 < x < \frac{1}{2}$,

$$f''(s) < 0, \qquad 1 < s < \frac{1}{2x}.$$

This proves (A.16).

(ii) The case $\frac{t}{x} \to \infty$ is obvious. For the case $\frac{t}{x} \to 1$, notice that $0 < x < t < 1$ and

$$\frac{(t-x)^3}{x^2} = o\left(\frac{(t-x)^2}{x}\right),$$

so

$$\begin{aligned} K^+(t,x) &= t \log \frac{t}{x} + (1-t) \log \frac{1-t}{1-x} \\ &= t\left[\left(\frac{t}{x} - 1\right) - \frac{1}{2}\left(\frac{t}{x} - 1\right)^2\right] + (1-t)\left[-\frac{t-x}{1-x} - \frac{1}{2}\left(\frac{t-x}{1-x}\right)^2\right] \\ &\quad + O\left(\frac{(t-x)^3}{x^2}\right) \\ &= \frac{1}{2}\frac{(t-x)^2}{x(1-x)}\left(1 + O\left(t + \frac{t}{x} - 1\right)\right). \end{aligned}$$ □

PROOF OF THEOREM 1.6. From (A.14) we have

$$\mathrm{BJ}_n^+ \leq \tfrac{1}{2}(\mathrm{HC}^*)^2,$$

so under $H_0$ the behavior of $\mathrm{BJ}_n^+$ is well controlled.

Now we consider the behavior of $\mathrm{BJ}_n^+$ under $H_1^{(n)}$. We examine the cases $r < \frac{\beta}{3}$ and $r > (1-\sqrt{1-\beta})^2$ separately; these two cases overlap and together cover the full region $\frac{1}{2} < \beta < 1$, $r > \rho^*(\beta)$.

First, for the case $r < \frac{\beta}{3}$ notice that $r < \frac{1}{4}$. Take $r_0$ such that $0 < r_0 < r$; as in Lemma A.4, it is easy to prove that, under $H_1^{(n)}$,

$$\max_{\{i:\, n^{-4r} < p_{(i)} < n^{-4r_0}\}} \left|\frac{p_{(i)}}{i/n} - 1\right| \to 0 \quad \text{in probability.}$$

Introduce the following statistic:

$$\mathrm{HC}^*_{r,r_0} = \max_{\{i:\, n^{-4r} < p_{(i)} < n^{-4r_0}\}} \left[\sqrt{n}\frac{i/n - p_{(i)}}{\sqrt{p_{(i)}(1-p_{(i)})}}\right].$$

Now, from (A.15)

$$n \cdot K^+(i/n, p_{(i)}) \sim \frac{1}{2}\left(\max\left\{\sqrt{n}\frac{i/n - p_{(i)}}{\sqrt{p_{(i)}(1-p_{(i)})}}, 0\right\}\right)^2 \quad \text{if } \frac{i/n}{p_{(i)}} \sim 1,$$

and so

$$(\mathrm{A}.17)\quad \mathrm{BJ}_n^+ \geq \max_{\{i:\, n^{-4r} < p_{(i)} < n^{-4r_0}\}} nK^+(i/n, p_{(i)}) = \tfrac{1}{2}[\mathrm{HC}^*_{r,r_0}]^2(1+o(1)).$$

Thus, in this case, $BJ_n^+$ is able to separate $H_0$ and $H_1^{(n)}$.

For the second case notice that $(r+\beta)/2\sqrt{r} < 1$. Pick a constant $q$ such that $\max\{(r+\beta)/2\sqrt{r}, \sqrt{r}\} < \sqrt{q} < 1$. Observe that under $H_1^{(n)}$,

$$\#\{i : p_i \leq n^{-q}\} \sim \mathrm{Binomial}(n, L_n n^{-[\beta+(\sqrt{q}-\sqrt{r})^2]}),$$
$$L_n n^{-[\beta+(\sqrt{q}-\sqrt{r})^2]} \gg n^{-q}.$$

This implies that, under $H_1^{(n)}$, for those $p$-values $p_{(i)} \sim n^{-q}$, $(i/n)/p_{(i)} \gg 1$; so, from (A.15),

$$nK^+(p_{(i)}, n^{-q}) \sim L_n \mathrm{Binomial}(n, L_n n^{-[\beta+(\sqrt{q}-\sqrt{r})^2]}).$$

As $1 - \beta - (\sqrt{q} - \sqrt{r})^2 > 0$, $BJ_n^+$ is able to separate $H_1^{(n)}$ and $H_0$. □

REMARK. For the second case $HC^*$ is more powerful than $BJ_n^+$. In fact, the best $BJ_n^+$ can do is to choose $i$ as large as possible while keeping $(i/n)/p_{(i)} \gg 1$. This is roughly equivalent to choosing $i \sim n^{1-q}$ with $q$ satisfying

$$n^{1-\beta} n^{-(\sqrt{q}-\sqrt{r})^2} \gg n^{1-q} \iff \sqrt{q} > (r+\beta)/(2\sqrt{r}).$$

As a result, $BJ_n^+ \approx L_n n^{1-(r+\beta)^2/(4r)}$. To see the main idea, take the region $\frac{1}{2} < \beta < \frac{3}{4}$, $\rho^*(\beta) < r < 1 - \beta$ for comparison. For $(r, \beta)$ in this range, $(HC^*)^2 \approx L_n n^{2r-\beta+1}$. Since $1 - (r+\beta)^2/(4r) < 2r - \beta + 1$, $HC^*$ has better performance than $BJ_n^+$.

**A.5. Proof of Lemma 5.1.** With $z_i \stackrel{\text{i.i.d.}}{\sim} N(0,1)$, $\delta_n = \mu_n^2 = 2r\log(n)$, $\chi^2(\delta_n) =_D z_1^2 + z_2^2 + \cdots + (z_\nu + \mu_n)^2$, so

$$P\{\chi_\nu^2(\delta_n) \geq 2q\log(n)\}$$
$$= P\{z_1^2 + z_2^2 + \cdots + (z_\nu + \mu_n)^2 \geq \sqrt{2q\log(n)}\}$$
$$= \frac{1}{(2\pi)^{\nu/2}} \int_{A(\theta_1,\rho)} \cos^{\nu-2}\theta_1 \rho^{\nu-1} e^{-\rho^2/2}\, d\theta_1\, d\rho$$
$$\times \int_{-\pi/2}^{\pi/2} \cos^{\nu-3}\theta_2\, d\theta_2 \cdots \int_0^{2\pi} d\theta_{\nu-1}$$
$$= \frac{2^{-\nu/2+1}}{\sqrt{\pi}\Gamma((\nu-1)/2)} \int_{|\theta_1|\leq\pi/2} \cos^{\nu-2}\theta_1\, d\theta_1$$
$$\times \int_{[-\sqrt{2r}\sin\theta_1 + \sqrt{2q-2r\cos^2\theta_1}]\sqrt{\log(n)}}^{\infty} \rho^{\nu-1} e^{-\rho^2/2}\, d\rho,$$

where $A(\theta_1, \rho) = \{(\theta_1, \rho) : |\theta_1| \leq \frac{\pi}{2}, \rho^2 + 2\rho\mu_n \sin\theta_1 \geq 2(q-r)\log(n)\}$.

The case $\nu = 1$ is obvious, while for $\nu = 2$ we have

$$P\{\chi_2^2(\delta_n) \geq 2q \log(n)\}$$
$$= \frac{1}{\pi} \int_{|\theta_1| \leq \pi/2} n^{-[\sqrt{q - r\cos^2\theta_1} - \sqrt{r}\sin\theta_1]^2} d\theta_1$$
$$= \frac{1}{\pi \sqrt{\log(n)}} n^{-[\sqrt{q} - \sqrt{r}]^2} \int_0^\infty e^{-\sqrt{r/q}(\sqrt{q} - \sqrt{r})^2 x^2} dx \, (1 + o(1))$$
$$= \frac{1}{\sqrt{2\pi \log(n)}} \frac{1}{\sqrt{2q} - \sqrt{2r}} \left(\frac{r}{q}\right)^{-1/4} \cdot n^{-[\sqrt{q} - \sqrt{r}]^2} (1 + o(1)).$$

Now consider the case $\nu \geq 3$. Notice that, for fixed $r$ and $q$ but large $n$, $n^{-[-\sqrt{r}\sin\theta_1 + \sqrt{q - r\cos^2\theta_1}]^2}$ obtains its maximum rate of growth at $\theta_1 = \frac{\pi}{2}$, and for $\theta_1 \approx \frac{\pi}{2}$,

$$[-\sqrt{r}\sin\theta_1 + \sqrt{q - r\cos^2\theta_1}]^2$$
$$\sim [\sqrt{q} - \sqrt{r}]^2 + \frac{\sqrt{r}}{\sqrt{q}}(\sqrt{q} - \sqrt{r})^2 \left(\theta_1 - \frac{\pi}{2}\right)^2;$$

moreover, notice that for any $y \to \infty$,

$$\int_y^\infty \rho^{\nu-1} e^{-\rho^2/2} d\rho = y^{\nu-2} e^{-y^2/2}(1 + o(1)).$$

These enable us to write

$$\int_{|\theta_1| \leq \pi/2} \cos^{\nu-2}\theta_1 \, d\theta_1 \int_{[-\sqrt{2r}\sin\theta_1 + \sqrt{2q - 2r\cos^2\theta_1}]\sqrt{\log(n)}}^\infty \rho^{\nu-1} e^{-\rho^2/2} d\rho$$
$$= ([\sqrt{2q} - \sqrt{2r}]\sqrt{\log(n)})^{\nu-2}$$
$$\quad \times \left[\int_{0 \leq \theta_1 \leq \pi/2} \cos^{\nu-2}\theta_1 \cdot n^{-[-\sqrt{r}\sin\theta_1 + \sqrt{q - r\cos^2\theta_1}]^2} d\theta_1\right] (1 + o(1))$$
$$= ([\sqrt{2q} - \sqrt{2r}]\sqrt{\log(n)})^{\nu-2}$$
$$\quad \times \left[\int_0^1 (1 - x^2)^{(\nu-3)/2} n^{-(\sqrt{q - r + rx^2} - \sqrt{r}x)^2} dx\right] (1 + o(1))$$
$$= ([\sqrt{2q} - \sqrt{2r}]\sqrt{\log(n)})^{\nu-2}$$
$$\quad \times \left[\int_0^1 (1 + x)^{(\nu-3)/2}(1 - x)^{(\nu-3)/2} n^{-(\sqrt{q - r + rx^2} - \sqrt{r}x)^2} dx\right] (1 + o(1))$$
$$= 2^{(\nu-3)/2}([\sqrt{2q} - \sqrt{2r}]\sqrt{\log(n)})^{\nu-2}$$
$$\quad \times \left[\int_0^1 (1 - x)^{(\nu-3)/2} n^{-(\sqrt{q - r + rx^2} - \sqrt{r}x)^2} dx\right] (1 + o(1)).$$

To evaluate this integration, notice that

$$\int_0^1 (1-x)^{(\nu-3)/2} n^{-(\sqrt{q-r+rx^2}-\sqrt{r}x)^2} \, dx$$

$$= \int_0^1 x^{(\nu-3)/2} n^{-(\sqrt{q-2rx+rx^2}-\sqrt{r}+\sqrt{r}x)^2} \, dx$$

$$= n^{-(\sqrt{q}-\sqrt{r})^2} \left[ \int_0^1 x^{(\nu-3)/2} n^{-2\sqrt{r/q}(\sqrt{q}-\sqrt{r})^2 x} \, dx \right] (1+o(1))$$

$$= \frac{1}{\log(n)} n^{-(\sqrt{q}-\sqrt{r})^2}$$

$$\times \left[ \int_0^{\log(n)} \left( \frac{x}{\log(n)} \right)^{(\nu-3)/2} e^{-\sqrt{r/q}(\sqrt{2q}-\sqrt{2r})^2 x} \, dx \right] (1+o(1))$$

$$= (\log(n))^{(1-\nu)/2} n^{-(\sqrt{q}-\sqrt{r})^2}$$

$$\times \left[ \int_0^\infty x^{(\nu-3)/2} e^{-\sqrt{r/q}(\sqrt{2q}-\sqrt{2r})^2 x} \, dx \right] (1+o(1))$$

$$= \Gamma\left(\frac{\nu-1}{2}\right) \left[ \sqrt{\frac{r}{q}}(\sqrt{2q}-\sqrt{2r})^2 \log(n) \right]^{(1-\nu)/2} n^{-(\sqrt{q}-\sqrt{r})^2}(1+o(1)).$$

Finally, we have

$$P\{\chi_\nu^2(\delta_n) \geq 2q\log(n)\}$$

$$= \frac{1}{\sqrt{2\pi \log(n)}} \left( \frac{r}{q} \right)^{(1-\nu)/4} \frac{1}{\sqrt{2q}-\sqrt{2r}} n^{-(\sqrt{q}-\sqrt{r})^2}(1+o(1)).$$

**Acknowledgments.** The authors thank Yoav Benjamini and Iain Johnstone for extensive discussions, references and encouragment. They also thank Howard Wainer, Betsy Becker, Henry Braun, Art Owen, David Siegmund and anonymous referees for useful references and pointers.

DEPARTMENT OF STATISTICS
STANFORD UNIVERSITY
STANFORD, CALIFORNIA 94305
USA
E-MAIL: donoho@stanford.edu